\documentclass[a4paper, 12pt]{article} 
\usepackage[utf8]{inputenc}				
\usepackage[numbers,square,sort&compress]{natbib}
\usepackage{epsfig,latexsym,amsfonts,amsmath,amssymb,color,mathrsfs}
\usepackage{setspace,dsfont}
\newtheorem{Theorem}{Theorem}
\newtheorem{Corollary}{Corollary}

\newtheorem{Lemma}{Lemma}

\newtheorem{Remark}{Remark}
\numberwithin{equation}{section}

\usepackage{amsfonts,amsmath,color}
\usepackage{graphicx}
\def\IntO{\int\limits_{\Omega}}           
\def\IntQ{\int\limits_{Q}}

\def\dvg{{\rm div}}
\def\Rd{{\mathbb R}^d}

\def\cP{{\mathcal P}}

\def\ep{{ \bf e}_\p}

\newcommand\be{\begin{eqnarray*}}    
\newcommand\ee{\end{eqnarray*}}      
\newcommand\ben{\begin{eqnarray}}
\newcommand\een{\end{eqnarray}}

\numberwithin{equation}{section}
\def\wt{\widetilde}

 \def\Normt#1{\mid\!\mid\!\mid #1 \mid\!\mid\!\mid}   

\def\e{{\rm e}}

\def\st{{1-2s}}
\def\y{{\mathbf y}}
\def\p{{\mathbf p}}

\def\cV{{\mathcal V}}

\begin{document}

\title{Functional a posteriori estimates\\ for the fractional Laplacian problem}
\author{Alexander Nazarov\footnote{ St. Petersburg Department of V.A. Steklov 
Institute of Mathematics,\\
\indent Fontanka 27, 191023, St. Petersburg, Russia; e-mail: 
nazarov@pdmi.ras.ru, repin@pdmi.ras.ru} ,
\setcounter{footnote}{0}
Sergey Repin\footnotemark
}

\maketitle

\hfill{\bf To the memory of Prof. R. Lazarov}

\abstract{
The paper is concerned with a posteriori estimates for approximations
of boundary value problems generated by the spectral fractional Laplace
operator. The derivation is based upon the  Stinga--Torrea extension, which
generalizes  the Caffarelli--Silvestre extension and
 transfers the corresponding nonlocal problem in a bounded domain to a local 
problem
of higher dimensionality. A posteriori estimates are first derived
for this local problem. Two--sided error bounds for the original problem 
follow from them.  The estimates are fully computable and contain
no conditions and constants depending on a method or mesh used to
compute an approximation. They are valid for any energy admissible 
approximation of the extended problem.

}
\maketitle
\large

\section{Introduction}
\label{intr}
Boundary value problems with fractional
differential operators generate new interesting problems in numerical analysis,
which require special methods (at this point we address the reader to
 \cite{BorNoc,ChDeSt,DuLaPa,HaLaMa,HaLaMaMa,Vab} and works cited therein).
After an approximation has been constructed, there rises the
question on its accuracy. 
A priori error estimates
give a theoretical idea of the accuracy that can be achieved provided that
all computations are performed exactly and a certain number of additional
assumptions are fulfilled. In the context of fractional Laplacian problem
they have been studied in \cite{Aco,BorNoc,No1} and some other publications.
Error estimates able to efficiently evaluate the quality of a particular
numerical solution are usually deduced within the framework of a
posteriori analysis.
 Currently,
there exist a wide range of a posteriori estimates from easily calculated error 
indicators
(commonly used in adaptive computations) 
up to fully guaranteed a posteriori estimates.
Among the estimates of the last mentioned type, an important class is formed by 
the so-called
{\em functional} a posteriori estimates.
They are derived with minimal assumptions on the solution and its approximation
and reflect the most general 
relations  that characterise  distances between the
(generalised) solution of the problem under investigation and
an arbitrary approximation in the corresponding energy class. 
Specific
properties of approximations or numerical methods and various regularity 
conditions are not used,
but this information can be utilised later when an a posteriori estimate of the 
functional
type is applied to a particular problem.
Theory of functional a posteriori estimates
is well studied for various problems generated by differential equations
of the local type (e.g., see
\cite{ReGruyter,ReSabook} ).

In this paper, we deduce such type estimates for the problem  $\cP^s_\Omega$.
Our considerations are based on the relations between problems
 $\cP^s_\Omega$ and  $\cP^s_Q$ established in \cite{StiTor} (for problems in 
the whole space see \cite{MO, CafSil}).
 They suggest an idea to reduce 
 analysis of the nonlocal problem 
$\cP^s_\Omega$ to
 analysis of local problem  $\cP^s_Q$ (which, however, has higher 
dimensionality).
 This problem can be used to generate approximations (e.g., see \cite{No1}).
 
 In Section \ref{apo}, we derive a posteriori relations able to control the 
accuracy
 of any approximation of the extended problem regardless of the method
 by which it has been constructed. Theorem~\ref{th1} establishes a posteriori 
error identity
 (\ref{IDmain})
 that serves as a basis for further analysis. Its left hand side is the sum of 
two error
 norms associated with approximations of the solution $w$ and the flux $\p$.
 The terms forming the right hand side are either known or can be estimated by 
known
 quantities and above mentioned norms. There are different ways to perform such 
an estimation.
 We apply the simplest method and obtain two sided estimates (\ref{id3a}) and 
(\ref{id4})
 that have fully computable right hand sides. If the problem $\cP^s_Q$ is solved
 by a mixed type method, then these estimates provide suitable error control 
tools.
Also, (\ref{IDmain}) has a reduced form (\ref{hyp}) with fully computable right 
hand
side. It arises if the last component of the flux $y_{d+1}$ is subject to the 
condition
(\ref{id3}).
Theorem \ref{th2} gives fully computable bounds for  energy norm of the error.

In Section \ref{ex}, the estimates are applied to spectral type approximations.
We show that in this case the condition (\ref{id3}) can be satisfied fairly 
easily
and deduce the corresponding error estimates (\ref{a7}) and (\ref{a8}).
They are verified in a series of numerical  tests (for the simplest case 
$s=\frac 12$).
Numerical tests show that the estimates 
adequately characterise the norms of deviations from the  solution 
for both exact and coarse approximations.
Note that the estimates used in the calculations
were obtained from error identities
using very simple transformations.
Using more sophisticated methods
will certainly lead to significantly more accurate estimates.
Also, it is worth noting that the method suggested in the paper can be extended
to other problems with fractional differential operators
provided that we have an extension of the corresponding nonlocal
problem in ${\mathbb R}^d$ to a local problem in ${\mathbb R}^{d+1}_+$ similar 
to (\ref{i2})--(\ref{i4}).

\section{Notation and preliminaries}
\label{not}
Throughout this paper,
 $\Omega$ is an open,  bounded and connected  domain in $\Rd$, $d\geq 1$
 with Lipschitz  boundary $\partial\Omega$ and  $0<s<1$.
We consider the following problem:
\medskip

{\bf Problem $\cP^s_\Omega$:} find
$u\in \widetilde H^s(\Omega)$ such that
\ben
\label{i1}
(-\Delta)^s_{\sf sp}u=f \quad\mbox{in}\quad \Omega.
\een
Here $f\in L^2(\Omega)$,
and $\widetilde H^s(\Omega)$ stands for the ``Dirichlet'' subspace of  the 
classical Sobolev--Slobodetskii space $H^s(\Rd)$ (e.g., see section 2.3.3 
of 
the book \cite{Tribel}), 
$$
\widetilde H^s(\Omega):=\{u\in H^s(\Rd)\ :\ u\equiv0 \quad\mbox{in}\quad 
\Rd\setminus\Omega\}.
$$
By $\|\cdot\|$ we denote norms of scalar and vector valued functions
in $L^2(\Omega)$ and   $(\cdot,\cdot)$ stands for the corresponding inner 
product.
Further, the {\em spectral fractional Laplacian} in (\ref{i1}) is defined by 
the relation
\be
(-\Delta)^s_{\sf sp}u=
\sum\limits^\infty_{j=1}\lambda^s_j(u,\phi_j)\phi_j,
\ee
where $\phi_j$ are the Dirichlet--Laplacian 
eigenfunctions (normalized in $L^2(\Omega)$) related to the
eigenvalues $\lambda_j$. In addition, 
we define the norm
$$
\|u\|^2_s:=\langle (-\Delta)^s_{\sf 
sp}u,u\rangle=\sum\limits^\infty_{ju=1}\lambda^s_j(u,\phi_j)^2,
$$
where $\langle\cdot,\cdot\rangle$ denotes the pairing between
$\widetilde H^s(\Omega)$ and its dual space $H^{-s}(\Omega)$.

One of the main difficulties is stipulated by non-locality of the operator
 $(-\Delta)^s_{\sf sp}$.
In \cite{CafSil} (see also \cite{MO}), it was shown that the fractional 
Laplacian problem can be 
studied as a
Dirichlet-to-Neumann operator by means of an extension problem 
in the half-space ${\mathbb R}^d\times(0,+\infty)$.
In \cite{StiTor}, the authors introduced similar extensions for a wider class
of problems including those defined in bounded domains. In the paper, we
use these results.

Let \,the \,added \,variable \,be \,denoted \,by $t$ and $w_t:=\frac{\partial 
w}{\partial 
t}$.
 Henceforth, \,($d+1$)-dimensional \,vectors \,are \,denoted \,by \,bold 
letters,
 e.g., $\y:=\{y_1,\dots,y_d,y_{d+1}\}$, 
 $\Delta_x$ stands for the standard Laplacian operator in ${\mathbb R}^d$ 
related
to the variables $x_1,x_2,\dots,x_d$, 
and
$$
\nabla_{xt}w:=\left\{\frac{\partial w}{\partial x_1},\dots,\frac{\partial 
w}{\partial x_d},
\frac{\partial w}{\partial t}\right\}=
\{\nabla_x w, w_t\}.
$$

The corresponding extended problem reads as follows.
\medskip

{\bf Problem $\cP^s_Q$:} find $w\in\mathcal V$ such that
\ben
\label{i2}
w_{tt}+\frac{1-2s}{t}w_t+\Delta_xw=0  &&\quad\mbox{in}\quad 
Q:=\Omega\times(0,+\infty),\\
\label{i3}-\lim\limits_{t\to0+}t^{1-2s}w_t=g  &&\quad\mbox{in}\quad \Omega,\\
\label{i4}
 w=0
&&\quad\mbox{on}\quad\Gamma:=\partial\Omega\times(0,+\infty).
\een
Here $g\in L^2(\Omega)$ is a given function,
and
$\mathcal V$ is the weighted Sobolev space in the half-cylinder $Q$ 
defined as the closure of smooth functions with bounded support
vanishing on $\Gamma$ with respect to the norm
\be
\|w\|^2_{\mathcal V}:=\Normt{\nabla_{xt}w}^2:=\IntQ 
t^\st|\nabla_{xt}w|^2\,dxdt.
\ee
We note that boundedness of this norm implies  that $\nabla_{xt} w$ must
decay with a certain rate as $t\rightarrow \infty$. In view of the homogeneous 
boundary
conditions on $\partial\Omega$, this means that $w$ must decay also.

It is convenient to introduce a ($d+1$)-dimensional vector valued function 
(flux)
$$\p:=t^{1-2s}\nabla_{xt}w.
$$
 Then,  (\ref{i2}) can be rewritten in the form
\be
\dvg_{xt}\p=0,\qquad 
\ee
and (\ref{i3}) is the Neumann condition 
\be
p_{d+1}(x,0)=-g(x)
\ee
stated on $\Omega$ (which is a part of the boundary
of $Q$).

 In what follows, we distinguish
scalar 
products in $\Rd$ and in $\mathbb R^{d+1}$ and denote them
by $\cdot$ and $\odot$, respectively (thus,  $|\y|^2:=\y\odot\y$)
and introduce the norms
\be
\Normt{w}^2:=\IntQ t^\st |w|^2 \,dxdt,\quad \Normt{\y}^2:=\IntQ t^\st |\y|^2 
\,dxdt
\ee
for the scalar and vector valued functions in $Q$. 

The solution of (\ref{i2})--(\ref{i4}) are understood in the standard weak sense 
by 
the relation
\ben
\label{gen}
\IntQ t^\st(w_t\eta_t+\nabla_xw\cdot\nabla_x\eta)\,dxdt-\IntO 
g\eta(x,0)\,dx=0,\qquad
\forall \eta\in {\mathcal V},
\een
or (in the terms of the flux)
\ben
\label{genmod}
\IntQ \p\odot\nabla_{xt}\eta\; \,dxdt-\IntO g\eta(x,0)\,dx=0,\qquad
\forall \eta\in {\mathcal V}.
\een

Connections between Problems $\cP^s_\Omega$ and $\cP^s_Q$ were thorough studied 
in \cite{StiTor}, where it was shown that  that if
\ben
\label{i6}
g=\,C_sf,\qquad\mbox{where}\qquad
C_s=2^{1-2s}\frac{\Gamma(1-s)}{\Gamma(s)},
\een
then the trace of the function $w$  equals the solution of 
$\cP^s_\Omega$, i.e.,
 $w(x,0)=u(x)$.\footnote{We outline that this correspondence takes place for 
$s\in (0,1)$ only and admits
 no continuation as $s\rightarrow 1$.}
 This fact opens a way to solve the problem numerically by means of methods
 developed for local type boundary value problems (e.g., see \cite{No1}).

 Below we deduce a posteriori estimates that can be applied for approximations
 obtained in this way.
\medskip

The problem (\ref{i2})--(\ref{i4}) has also a variational statement.
Let us define
$$
J(\wt w):=\frac12\Normt{\nabla_{xt}\wt w}^2-\IntO g\wt w(x,0) \,dx.
$$ 
and consider the 
variational problem: find $w\in {\mathcal V}$ such that
\ben
\label{i7}
J(w)=\inf\limits_{{\wt w}\in {\mathcal V}} J({\wt w}).
\een
Existence of the minimizer $w$ follows from standard  convexity and coercivity 
arguments.
It is easy to see that (\ref{gen}) is just the Euler equation for the problem 
(\ref{i7}).

\begin{Remark}
By the same reason, any $w\in\mathcal V$ solving the equation (\ref{i2}) 
satisfies
\ben
\label{i8}
\Normt{\nabla_{xt}w}^2=\inf\limits_{{\wt w}\in {\mathcal V}
\atop \wt w(\cdot,0)= w(\cdot,0)} \Normt{\nabla_{xt}\wt w}^2.
\een
\end{Remark}

\begin{Lemma}
For any $\wt w \in {\mathcal V}$ it holds
\ben
\label{IDM}
\frac12\Normt{\nabla_{xt}(\wt w -w)}^2=J(\wt w )-J(w).
\een
\end{Lemma}

{\bf Proof}.
We have
\begin{multline*}
J(\wt w )-J(w)=\frac12\Normt{\nabla_{xt}\wt w }^2-\frac12\Normt{\nabla_{xt} 
w}^2-
\IntO g(\wt w(x,0) -w(x,0)) \,dx\\
=\frac12\Normt{\nabla_{xt}(\wt w -w)}^2+
\IntQ t^\st\nabla_{xt}w\odot\nabla_{xt}(\wt w -w)\,dxdt\\
-\IntO g(\wt w(x,0) 
-w(x,0))\,dx.
\end{multline*}
Setting $\eta=\wt w-w$ in (\ref{gen}) we
arrive at (\ref{IDM}).
\hfill$\square$\medskip

\begin{Lemma}
For any $\wt w \in {\mathcal V}$ it holds
\ben
\label{ID2}
C_s\|\wt w(\cdot,0)\|_s^2\le\, \Normt{\nabla_{xt}\wt w}^2.
\een
and
\ben
\label{ID1}
\Normt{\wt w}\,\leq\, 
C_F 
\Normt{\nabla_{xt}\wt w},
\een
where $C^2_F(\Omega)=\lambda_1^{-1}$ and $\lambda_1$ is the minimal eigenvalue
of the Dirichlet Laplacian in $\Omega$.

Moreover, for any $w \in {\mathcal V}$ solving the equation 
(\ref{i2}) the inequality (\ref{ID2}) 
holds as the equality, i.e.,
\ben
\label{ID3}
C_s\|w(\cdot,0)\|_s^2=\,\Normt{\nabla_{xt}w}^2.
\een
\end{Lemma}
{\bf Proof}.
Let $\wt w\in\mathcal V$. From the classical Steklov--Friedrichs inequality we 
have
\be
\|\wt w(\cdot,t)\|^2\leq C^2_F(\Omega)\|\nabla_x \wt w(\cdot,t)\|^2\qquad 
\forall 
\;t\in (0,\infty).
\ee
Hence we obtain
 a weighted Friedrichs type inequality
\be
\IntQ t^\st|\wt w(x,t)|^2\,dxdt\leq C^2_F(\Omega)\IntQ t^\st|\nabla_{xt} 
\wt w(x,t)|^2
\,dxdt,
\ee
which yields (\ref{ID1}).

Next, let $w$ be the solution of (\ref{i2})--(\ref{i4}) for some $g$. Define 
$f$ in 
accordance with (\ref{i6}) and recall that $u=w(\cdot,0)$ solves (\ref{i1}). 
Substituting $\wt w=w$ in (\ref{gen}) we get
$$
\aligned
\Normt{\nabla_{xt}w}^2= &\IntQ t^\st(|w_t|^2+|\nabla_xw|^2)\,dxdt=\IntO g 
w(x,0)\,dx\\
= &\, C_s\IntO f w(x,0)\,dx
= C_s\IntO f u \,dx=C_s\IntO ((-\Delta)^s_{\sf sp}u) u \,dx\\ 
= &\, C_s\|u\|^2_s=C_s\|w(\cdot,0)\|^2_s,
\endaligned
$$
and (\ref{ID3}) follows.

Finally, we notice that (\ref{ID2}) follows from (\ref{i8}) and (\ref{ID3}).
\hfill$\square$\medskip

\begin{Corollary}
Since for every $v\in \widetilde H^s(\Omega)$
$$
\|v\|^2:=
\sum\limits^\infty_{j=1}(v,\phi_j)^2\leq \frac{1}{\lambda^s_1}
\sum\limits^\infty_{k=1}\lambda^s_j(v,\phi_j)^2=C^{2s}_F(\Omega)\|
v\|^2_s,
$$
we conclude from (\ref{ID2}) that for any $\wt w\in\cV$
\ben
\label{ID4}
\|\wt w(\cdot,0)\|\leq C^s_F(\Omega)\|\wt w(\cdot,0)\|_s\leq 
C^s_F(\Omega)\kappa_s\Normt{\nabla_{xt}\wt w},\quad
\kappa_s=C^{-\frac 12}_s.
\een
\end{Corollary}

\begin{figure}
\begin{center}
{\includegraphics[height=6cm,width=12cm]{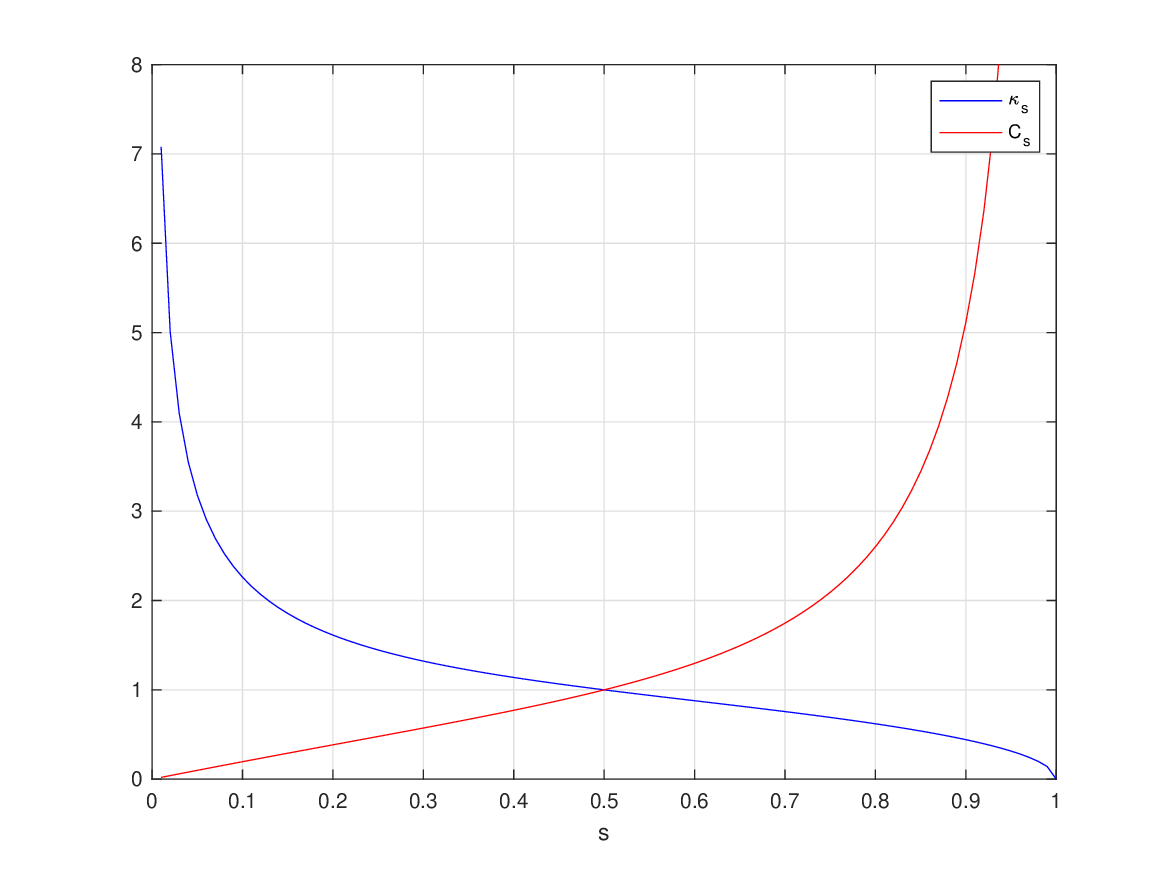}}
\end{center}
\caption{Constants $C_s$ and $\kappa_s$}
\label{figcon}
\end{figure}

The graphs of $C_s$ and $\kappa_s$ are given at Fig.~\ref{figcon}.

\section{A posteriori estimates}
\label{apo}
\subsection{A posteriori error identity for the problem  $\cP^s_Q$}
For the majority of boundary value problems errors are usually
evaluated in terms of combined norms (measures) that account
approximations of the primal and dual variables (e.g., approximations of the 
solution
function and of its flux).  For the divergence type equations the corresponding 
estimates follow from
error identities that hold for any approximation in the
energy space (see \cite{ReGruyter,Re2021}). In our case, this space is $\cV$. 
The
 space $Y^s$ of admissible 
fluxes consists of vector valued functions $\y$
such that  $\Normt{\y}<\infty$. In addition, we define
its subspace
\be
Y^s_{\dvg}:=
\Bigl\{\y\in Y^s\,\mid\,
\IntQ t^{2s-1}|\dvg_{xt}\y|^2\,dxdt<+\infty\Bigr\}.
\ee
Theorem below establishes the  error identity that holds for
approximations of  the Problem $\cP^s_Q$.

\begin{Theorem}
\label{th1}
For any ${\wt w}\in {\mathcal V}$ and any $\y\in Y^s_{\dvg}$
it holds
\begin{multline}
\label{IDmain}
\Normt{\nabla_{xt}{\e_w}}^2+\Normt{t^{2s-1}\ep}^2\,=\,
\Normt{\nabla_{xt}{\wt w}-t^{2s-1}\y}^2\\
+2\IntQ {\e_w} \dvg_{xt}\y  \,dxdt
-2\IntO {\e_w}(x,0)(g+y_{d+1}(x,0))\,dx,
\end{multline}
where $\e_w:=\wt w-w$ and $\ep:=\y-\p$.
\end{Theorem}
{\bf Proof}.
We have
\begin{multline}
\label{id1}
\Normt{\nabla_{xt}{\e_w}}^2+\Normt{t^{2s-1}\ep}^2\,=\,\Normt{\nabla_{xt}{\e_w}-
t^{2s-1}\ep}^2+\,2\IntQ \nabla_{xt} {\e_w}\odot\ep \,dxdt\\
=\,
\Normt{\nabla_{xt}{\wt w}-t^{2s-1}\y}^2+\,2\IntQ \nabla_{xt} {\e_w}\odot 
(\y-\p) \,dxdt.
\end{multline}
In view of (\ref{genmod}),
\begin{multline}
\label{id2}
\IntQ \nabla_{xt} {\e_w}\odot (\y-\p) \,dxdt
\\
=- \IntQ {\e_w}\dvg_{xt}\y  \,dxdt+
\IntO {\e_w}(x,0)(p_{d+1}-y_{d+1}(x,0))\,dx\\
=\IntQ {\e_w}\dvg_{xt}\y  \,dxdt-
\IntO {\e_w}(x,0)(g+y_{d+1}(x,0))\,dx,
\end{multline}
and (\ref{IDmain}) follows from (\ref{id1}) and (\ref{id2}).
\hfill$\square$\medskip

The identity (\ref{IDmain}) has a reduced form if the function $\y$ is selected 
in a special way. Namely, let
\be
Y^s_g:=\big\{\y\in \; Y^s\,\mid\,\dvg_{xt}\y=0 \ \ \mbox{in}\ Q,\quad 
y_{d+1}(x,0)=-g(x) \ \ \mbox{in}\ \Omega
\big\}.
\ee
Easily, if we denote $ \y_d:=\{y_1,y_2,\dots,y_d\}$, then this definition can 
be 
rewritten as follows:
\ben
\label{id3}
Y^s_g=\Big\{\y\in\;Y^s\,\mid\,y_{d+1}(x,t)=-g(x)-\int\limits^t_0\dvg_x \y_d \,
dt\Big\}.
\een

\begin{Corollary}
For any ${\wt w}\in {\mathcal V}$ and any $\y\in Y^s_g$ it holds
\ben
\label{hyp}
\Normt{\nabla_{xt}{\e_w}}^2+\Normt{t^{2s-1}\ep}^2=
\Normt{\nabla_{xt}{\wt w}-t^{2s-1}\y}^2.
\een
\end{Corollary}

Indeed,
if $\y\in Y^s_g$, then both 
integrals in the right hand side of (\ref{IDmain}) vanish,
and we arrive at the reduced identity (\ref{hyp}). 
From (\ref{hyp}) it follows that
\ben
\label{H1}
\Normt{\nabla_{xt}{\e_w}}^2\;\leq\;
\Normt{\nabla_{xt}{\wt w}-t^{2s-1}\y}^2\qquad \forall \y\in Y^s_g.
\een
This inequality
can be viewed as a hypercircle type error estimate
for the problem $\cP^s_Q$.
Notice that $\p\in Y^s_g$, so that
\be
\Normt{\nabla_{xt}{\e_w}}^2\;=\;\inf\limits_{\y\in Y^s_g}
\Normt{\nabla_{xt}{\wt w}-t^{2s-1}\y}^2.
\ee

\subsection{Computable estimates of the norm $\Normt{\nabla_{xt}\e_w}$}

\indent
Right hand sides of (\ref{hyp}) and (\ref{H1}) are fully 
computable. However, these
relations hold
only for  $\y$ whose last component is defined in a special way (\ref{id3}). 
In essence, this condition comes down to calculating the primitive with respect 
to time for a given function. 
It does not seem difficult to fulfil.
However,  the identity (\ref{IDmain}) implies more general estimates that
operate with  $\y\in Y^s_{\dvg}$. Its right hand side can be estimated
from above by the norm $\Normt{\nabla_{xt}\e_w}$. Therefore,
two--sided a posteriori estimates for the combined (primal--dual) error norm 
follow
from  (\ref{IDmain}).

There are different methods to estimate the integrals in the right hand side of
(\ref{IDmain}).
Below we use the 
simplest one, which yields simple error bounds, which may
overestimate the error.
We believe that more sophisticated methods (as those discussed in
\cite{Re2021} for local type problems of elliptic type) will generate much 
sharper
estimates. However, for now the goal is to show principal ability to deduce
fully computable error estimates for any $\wt w\in \cV$ and $\y\in Y^s_{\dvg}$.
Therefore, in this paper we confine ourselves to the simplest variant.
 
In view of (\ref{ID1}), for any positive $\alpha_1$ we have 
\begin{multline*}
2\Big|\IntQ {\e_w} \dvg_{xt}\y \,dxdt\Big|\leq 
\frac{\alpha_1}{C^2_F(\Omega)}\Normt{\e_w}^2+
\frac{C^2_F(\Omega)}{\alpha_1}\Normt{ t^{2s-1}\dvg_{xt}\y}^2
\\
\leq
\alpha_1\Normt{{\nabla_{xt}\e_w}}^2+
\frac{C^2_F(\Omega)}{\alpha_1}\Normt{ t^{2s-1}\dvg_{xt}\y}^2.
\end{multline*}
Similarly, by (\ref{ID4}) we have  for any positive $\alpha_2$
\begin{multline*}
2\Big|\IntO {\e_w}(x,0)(g+y_{d+1}(x,0))\,dx\Big|\\
\leq  
\frac{\alpha_2}{C^{2s}_F(\Omega)\kappa_s^2}
\,\|{\e_w}(\cdot,0)\|^2+
\frac{C^{2s}_F(\Omega)\kappa_s^2}{ \alpha_2}\,\|g+y_{d+1}(x,0)\|^2
\\
\leq
{\alpha_2}\Normt{\nabla_{xt} {\e_w}}^2+
\frac{C^{2s}_F(\Omega)\kappa_s^2}{ \alpha_2}\,\|g+y_{d+1}(x,0)\|^2.
\end{multline*}
Hence if
$\alpha_1+\alpha_2< 1$  then (\ref{IDmain}) implies the 
estimate
\begin{multline}
\label{id3a}
(1-\alpha_1-\alpha_2)\Normt{\nabla_{xt}{\e_w}}^2+
\Normt{t^{2s-1}\ep}^2\\
\leq
\,\Normt{\nabla_{xt}{\wt w}-t^{2s-1}\y}^2
+\frac{C^2_F(\Omega)}{\alpha_1}\Normt{ t^{2s-1}\dvg_{xt}\y}^2\\
+\frac{C^{2s}_F(\Omega)\kappa_s^2}{ \alpha_2}\|g+y_{d+1}(\cdot,0)\|^2.
\end{multline}
Analogously we obtain the estimate from below
\begin{multline}
\label{id4}
(1+\alpha_1+\alpha_2)\Normt{\nabla_{xt}{\e_w}}^2+\|\ep
\|^2
\geq\, 
\Normt{\nabla_{xt}{\wt w}-t^{2s-1}\y}^2 \\
-\frac{C^2_F(\Omega)}{\alpha_1}\Normt{ t^{2s-1}\dvg_{xt}\y}^2
-\frac{C^{2s}_F(\Omega)\kappa_s^2}{ \alpha_2}\|g+y_{d+1}(\cdot,0)\|^2,
\end{multline}
which is useful for $\alpha_1$ and $\alpha_2$ large enough.

Estimates (\ref{id3a}) and (\ref{id4}) hold for the combined primal--dual error 
norm.
 Now, our goal is to deduce computable bounds of the norm 
$\Normt{\nabla_{xt}\e_w}$.


\begin{Theorem}
\label{th2}
For any $\y\in Y^s_\dvg$ and  any $\eta\in {\mathcal V}$ 
it holds
\ben
\label{maj1}
M_\ominus({\wt w};\eta)\;\leq\;
\Normt{\nabla_{xt}{\e_w}}\;\leq\;M_\oplus({\wt w};\y),
\een
where
\begin{multline*}
M_\oplus({\wt w};\y):=\,
\Normt{\nabla_{xt}{\wt w}-t^{2s-1}\y}\\
+\,
C_F 
\Normt{t^{2s-1}\dvg_{xt}\y}
+\,C^s_F\kappa_s\|y_{d+1}(\cdot,0)+g\|
\end{multline*}
and
\be
M_\ominus^2({\wt w};\eta):=
2\IntQ t^{1-2s}\nabla_{xt}{\wt w}\odot\nabla_{xt}\eta \,dxdt-2\IntO g\eta 
\,dx-\Normt{\nabla_{xt} \eta}^2.
\ee
Moreover,
\ben
\label{maj}
\sup\limits_{\eta\in {\mathcal V}}\;M_\ominus({\wt w};\eta)\;=\;
\Normt{\nabla_{xt}{\e_w}}\;=\;\inf\limits_{\y\in Y^s_{\dvg}}M_\oplus({\wt 
w};\y).
\een
\end{Theorem}
{\bf Proof}.
Consider the right part of (\ref{maj}) first.
From (\ref{gen}) we have for any $\eta\in {\mathcal V}$
\begin{multline}
\label{e1}
\IntQ t^\st\;\nabla_{xt} {\e_w}\odot\nabla_{xt}\eta\,dxdt \\
=\IntQ t^\st\;
\nabla_{xt}{\wt w}\odot\nabla_{xt} \eta\,dxdt-
\IntO g\eta(x,0)\,dx.
\end{multline}
Setting  $\eta={\e_w}$ and noticing that
\begin{multline*}
\IntQ 
({\e_w}\dvg_{xt}\y+\nabla_{xt}{\e_w}\odot\y)\,dxdt=\IntQ 
\dvg_{xt}({\e_w}\y)\,dxdt \\
=
- \IntO  {\e_w}(x,0)y_{d+1}(x,0)\,dx, 
\end{multline*}
we rearrange (\ref{e1})  as follows
\begin{multline}
\label{e2}
\Normt{\nabla_{xt}{\e_w}}^2=
\IntQ t^\st\nabla_{xt}{\wt w}\odot\nabla_{xt} {\e_w}\;\,dxdt-
\IntO g {\e_w}(x,0)\,dx\\
=\IntQ (t^\st\nabla_{xt}{\wt w}-\y)\odot\nabla_{xt} {\e_w}\;\,dxdt-\IntQ \e_w
\dvg_{xt}\y \,dxdt\\
-\IntO  {\e_w}(x,0)(y_{d+1}(x,0)+g)\,dx.
\end{multline}
Consider the first integral in the right hand side of (\ref{e2}). We have
\begin{multline*}
\Big|\IntQ (t^\st\nabla_{xt}{\wt w}-\y)\odot\nabla_{xt} {\e_w}\;\,dxdt\Big| \\
=
\Big|\IntQ t^\st(\nabla_{xt}{\wt w}-t^{2s-1}\y)\odot\nabla_{xt} 
{\e_w}\;\,dxdt\Big|
\leq\,
\Normt{\nabla_{xt}{\wt w}-t^{2s-1}\y}
\ \Normt{\nabla_{xt}{\e_w}}.
\end{multline*}
Next, by the Cauchy--Bunyakovski inequality and  (\ref{ID1}), we have
\begin{multline*}
\Big|\IntQ {\e_w}\dvg_{xt}\y \,dxdt\Big|\leq \,
\Normt{{\e_w}}\left(\IntQ t^{2s-1}|\dvg_{xt}\y|^2\,dxdt\right)^{\frac 12}
\\
\leq 
C_F\Normt{\nabla_{xt}{\e_w}}\ \Normt{t^{2s-1}\dvg_{xt}\y},
\end{multline*}
and by (\ref{ID4}),
\begin{multline*}
\Big|\IntO  {\e_w}(x,0)(y_{d+1}(x,0)+g)\,dx\Big|\leq 
\|{\e_w}(\cdot,0)\|\ \|y_{d+1}(\cdot,0)+g\| 
\\
\le C^s_F\kappa_s\Normt{\nabla_{xt}{\e_w}}\,\|y_{d+1}(\cdot,0)+g\|.
\end{multline*}
Hence we arrive at the estimate
\begin{multline}
\label{e8}
\Normt{\nabla_{xt}{\e_w}}\;\leq\;\Normt{\nabla_{xt}{\wt w}-t^{2s-1}\y} \\
+\,C_F\Normt{t^{2s-1}\dvg_{xt}\y}+\,C^s_F\kappa_s\|y_{d+1}(\cdot,0)+g\|,
\end{multline}
which is just the right inequality in (\ref{maj1}).

Set $\y=\p$ in (\ref{e8}). Then the last two terms vanish and 
$t^{2s-1}\y=\nabla_{xt}w$.
In this case, (\ref{e8}) holds as the equality. Thus, the
right equality in (\ref{maj}) is justified.
\medskip

Let us now prove the left equality.
Taking into account (\ref{IDM}) and the fact that $w$ is the minimizer of the 
functional $J$, we see that for any $\eta\in {\mathcal V}$
\begin{multline}
\label{IDMa}
\frac12\Normt{\nabla_{xt}({\wt w}-w)}^2=J({\wt w})-J(w)\geq
J({\wt w})-J({\wt w}-\eta)\\
=
\IntQ t^{1-2s}\nabla_{xt}{\wt w}\odot\nabla_{xt}\eta 
\,dxdt-\frac12\Normt{\nabla_{xt} \eta}^2-
\IntO g\eta \,dx=\frac12 M_\ominus^2(\wt w;\eta).
\end{multline}
On the other hand, if we set $\eta={\wt w}-w$ in (\ref{IDMa}) it becomes 
equality. This justifies the left equality in (\ref{maj}). 
\hfill$\square$\medskip
\begin{Remark}
Assume that $M_\oplus(\wt w;\y)=0$.
Since all terms  of the majorant vanish, we conclude that
\be
\y=t^\st\nabla_{xt}{\wt w},\quad \dvg_{xt}\y=0,\quad y_{d+1}=t^\st {\wt w_t}=-g,
\ee
what means that $\wt w=w$ and $\y=\p$.
\end{Remark}


\begin{Remark}
By means of (\ref{maj1}), we can construct a sequence of computable two sided
bounds of the error
using finite dimensional subspaces ${\mathcal V}_n\subset {\mathcal V}$  (${\rm 
dim} {\mathcal V}_n=n$)
and $Y_k\subset Y^s$ (${\rm dim} Y_k=k$).
If we denote by $\eta_n\in{\mathcal V}_n$ and $\y_k\in Y_k$ corresponding 
extremal elements, we have
\begin{multline*}
M_\ominus({\wt w};\eta_n)=
\sup\limits_{\eta\in {\mathcal V}_n}\;M_\ominus(\wt 
w;\eta)\leq\;\Normt{\nabla_{xt}({\wt w}-w)}\\
\leq
M_\oplus(\wt w;\y_k)=\inf\limits_{\y\in Y_k}M_\oplus(\wt w;\y).
\end{multline*}
It is clear that $M_\ominus({\wt w};\eta_{n+1})\geq M_\ominus({\wt w},\eta_n)$ 
if ${\mathcal 
V}_n\subset {\mathcal V}_{n+1}$,
and 
$$
\lim\limits_{n\rightarrow+\infty}M_\ominus({\wt w};\eta_n)
\;=\; \Normt{\nabla_{xt}({\wt w}-w)} 
$$
 if the spaces $\mathcal V_n$ are limit dense in ${\mathcal V}$.
 Analogously, 
 $$
 \lim\limits_{n\rightarrow+\infty}M_\oplus({\wt w};\y_k)
\;=\; \Normt{\nabla_{xt}({\wt w}-w)} 
$$
 if the spaces $Y_k$ are nested and limit dense in $Y^s_{\dvg}$. Hence 
theoretically 
the
 error can be determined with any desired accuracy. Certainly, in practice 
 admissible expenditures are always limited
 so that above presented error bounds may be not sharp.
  \end{Remark}
 
\subsection{A posteriori error estimates for the problem $\cP^s_\Omega$}
Above derived a posteriori estimates for approximations of the extended problem
 $\cP^s_Q$ open a way to majorate the difference
 between ${v}(x):=\wt w(x,0)$ and $u(x)$.
Since
${v}(x)-u(x)=\wt w(x,0)-w(x,0)$,
computable bounds of  
$$
\e_{u}:=v-u
$$
follow from the above derived estimates.
By (\ref{ID2}) we have
\ben
\label{e6a}
\|{\e_u}\|_{s,\Omega}\leq \kappa_s\Normt{\nabla_{xt} \e_w}.
\een
Hence (\ref{maj1}) implies the estimate
\ben
\label{majful}
\|{\e_u}\|_{s,\Omega}\,\leq \;\kappa_s\,M_\oplus({\wt w};\y).
\een
The right hand side of (\ref{majful}) contains known
functions ${\wt w}$ (approximate solution of $\cP_Q$) and
$\y$ (approximate flux of $\cP_Q$).

If $\wt w$ satisfies equation (\ref{i2}), then due to (\ref{ID3})  we
have equality in (\ref{e6a}). By (\ref{maj1}) we obtain
\be
\|{\e_u}\|_{s,\Omega}\,\geq \;\kappa_s\,M_\ominus({\wt w},\eta).
\ee

\begin{Remark}
Serious computational difficulties should be expected as we approach
both limiting cases
 $s= 0$ and $s= 1$. 
They are caused by the singularity of the coefficients in (\ref{i2}) and 
(\ref{i3}).
 The problem $\cP_Q$ degenerates at $s\rightarrow 1$
because the weight $t^{1-2s}$ becomes locally unsummable for $s=1$.
In addition, if $s\rightarrow 0$, then  the constant $\kappa_s$ blows up (see 
Fig. 
\ref{figcon}) so that  efficiency of the estimate (\ref{e6a}) deteriorates if 
$s$ approaches
$0$.
\end{Remark}

\section{Example}
\label{ex}

Consider our  problem with $s=\frac12$. Then $C_s=\kappa_s=1$. 
Let the 
solution of  $\cP^s_\Omega$
 be presented in the  form $u(x)=\sum\limits^\infty_{j=1} a_j \phi_j(x)$ 
(recall that
$\phi_j$ is the normalized Dirichlet--Laplacian eigenfunction associated
with the eigenvalue $\lambda_j$).
Since
\ben
\label{a1}
f=(-\Delta)_{\sf sp}^{\frac12} 
u=\sum\limits^\infty_{j=1}\lambda^{\,\frac 12}_j(u,\phi_j)\phi_j
=\sum\limits^\infty_{j=1} \lambda^{\,\frac 12}_j a_j\phi_j
\een
we conclude that $a_j=\lambda^{-\frac 12}_j\zeta_j$, where 
$$\zeta_j:=(f,\phi_j)=\IntO f\phi_j \,dx.
$$

Find $w$ (solution of  the problem $\cP^s_Q$) in the form
$w(x,t)=\sum\limits^\infty_{j=1} \beta_j(t)\phi_j(x)$.
It meets the relations
\ben
\label{a2}
w_{tt}+\Delta_x w=0\qquad\text{and}\qquad 
w_t(\cdot,0)=-g=-f.
\een

The differential equation in (\ref{a2}) and the assumption $w\in \mathcal V$ 
imply
$\beta_j(t)=c_j e^{-\lambda^{\,\frac 12}_jt}$.
Hence
\be
w_t(x,t)=-\sum\limits^\infty_{j=1} c_j \lambda^{\,\frac 12}_j 
e^{-\lambda^{\,\frac 12}_jt}\phi_j(x)
\ \
\text{and}\ \
w_t(x,0)=-\sum\limits^\infty_{j=1} c_j \lambda^{\,\frac 12}_j \phi_j(x)=-f.
\ee
Comparing the last relation with (\ref{a1}), we see that $c_j=a_j$. Therefore,
\ben
\label{a3}
w(x,t)=\sum\limits^\infty_{j=1}\lambda^{-\frac 12}_j\zeta_j \phi_j(x) 
e^{-\lambda^{\,\frac 12}_jt}
\een
and
\be
\p=\{\nabla_x w,w_t\}=
\left\{
 \sum\limits^\infty_{j=1} 
 \lambda^{-\frac 12}_j\zeta_j \nabla_x\phi_j(x) 
 e^{-\lambda^{\,\frac 12}_jt},
-\sum\limits^\infty_{j=1} 
\zeta_j \phi_j(x) 
e^{-\lambda^{\,\frac 12}_jt}
\right\}.
\ee

Thus, if $\phi_j$ and $\lambda_j$ are known, then theoretically
 solutions of both problems
$\cP_\Omega$ and $\cP_Q$ can be presented as series.
However,  $\phi_j$ and $\lambda_j$ are  known exactly for a limited amount
of special domains $\Omega$ so that in general
instead of them we have certain approximations $\psi_j$ and $\theta_j$, 
respectively.
Besides, we are forced to replace infinite series be finite sums.
Hence instead of the exact solution we have a spectral type approximation 
 in the form similar to (\ref{a3})
\ben
\label{a4}
\wt w(x,t)=\sum\limits
^N
_{j=1}\theta^{-\frac 12}_j\gamma_j \psi_j(x)  e^{-\theta^{\,\frac 12}_jt},\qquad
\gamma_j=(f,\psi_j),
\een
with certain positive nondecreasing numbers $\theta_j$  and functions $\psi_j$ 
vanishing on $\partial\Omega$. 

In this case, it is natural to find $\y=\{\y_d,y_{d+1}\}$ in the form
\be
\y_d=\sum\limits^N_{j=1}\Upsilon_j(x) e^{-\theta^{\,\frac 12}_jt},\qquad
y_{d+1}=\sum\limits^N_{j=1}{\theta^{-\frac 12}_j}\dvg_x\Upsilon_j(x) 
e^{-\theta^{\,\frac 12}_jt}.
\ee
Here $\Upsilon_j$
are certain vector valued functions such that $\dvg_x\Upsilon_j\in L^2(\Omega)$.
It is easy to see that so defined vector function $\y$ satisfies
the relation $\dvg_{xt}\y=0$. 

We have
\be
&&\nabla_x\wt w-\y_d=\sum\limits
^N
_{j=1}
 Q_j(x) e^{-\theta^{\,\frac 12}_j t},\qquad
Q_j(x):=\theta^{-\frac 12}_j
\gamma_j\nabla\psi_j(x)-\Upsilon_j(x),\\
&&\wt w_t-y_{d+1}=-\sum\limits^N_{j=1}
R_j(x)e^{-\theta^{\,\frac 12}_jt},\quad
R_j(x):=\gamma_j\psi_j(x)+\theta^{-\frac 12}_j\dvg_x \Upsilon_j(x).
\ee
Hence
\begin{multline*}
\Normt{\nabla_{xt}{\wt w}-\y}^2
=\sum\limits^N_{j=1}\sum\limits^N_{k=1}
\IntO (Q_j\cdot Q_k+R_jR_k)\,dx\int\limits^\infty_0  
 e^{-(\theta^{\,\frac 12}_j+\theta^{\,\frac 12}_k)t}dt\\
 =
 \sum\limits^N_{j=1}\sum\limits^N_{k=1}
\frac{1}{\theta^{\,\frac 12}_j+\theta^{\,\frac 12}_k}\IntO (Q_j\cdot 
Q_k+R_jR_k)\,dx
=:S^2_N,
\end{multline*}
and  the majorant in (\ref{maj1}) has the form
\be
M_\oplus(\wt w;\y)=S_N+
C^{\,\frac 12}_F\left(\,
\IntO\Big|\sum\limits^N_{j=1}{\theta^{-\frac 12}_j}\dvg_x\Upsilon_j+f\Big|^2\,dx
\right)^{\frac 12}.
\ee
In particular, if
we set $\Upsilon_j=\gamma_j\theta^{-\frac 12}_j\nabla\psi_j(x)$, then 
$Q_j(x)=0$ 
and
$$
R_j(x)=\gamma_j\theta^{-1}_j\rho_j(x), \qquad 
\rho_j(x):=\Delta_x\psi_j(x)+\theta_j\psi_j(x).
$$
In this case,
\ben
\label{SN}
S^2_N=\sum\limits^N_{j=1}\sum\limits^N_{k=1}
\frac{\gamma_j\gamma_k}
{{\theta_j}{\theta_k}(\theta^{\,\frac 12}_j+\theta^{\,\frac 12}_k)}
\IntO\rho_j\rho_k \,dx
\een
and
\ben
\label{a7}
M_\oplus(\wt w;\y)=S_N
+C^{\,\frac 12}_F
\left(\,\IntO\Big|\sum\limits^N_{j=1}\gamma_j{\theta^{-1}_j}\Delta_x\psi_j+f\Big
|^2\, dx
\right)^{\frac 12}.
\een
Also, we can use the estimate (\ref{id3a}). For the above defined $\y$
it reads
\ben
\label{a8}
(1-\alpha)\Normt{\nabla_{xt}\e_w}^2+\Normt{\ep}^2\leq S^2_N+
\frac{C_F}{\alpha} 
\IntO\Big|\sum\limits^N_{j=1}{\theta^{-1}_j}\gamma_j\Delta_x\psi_j+f\Big|^2\,dx.
\een

\begin{figure}
\begin{flushleft}
\includegraphics[height=8.9cm,width=16cm]{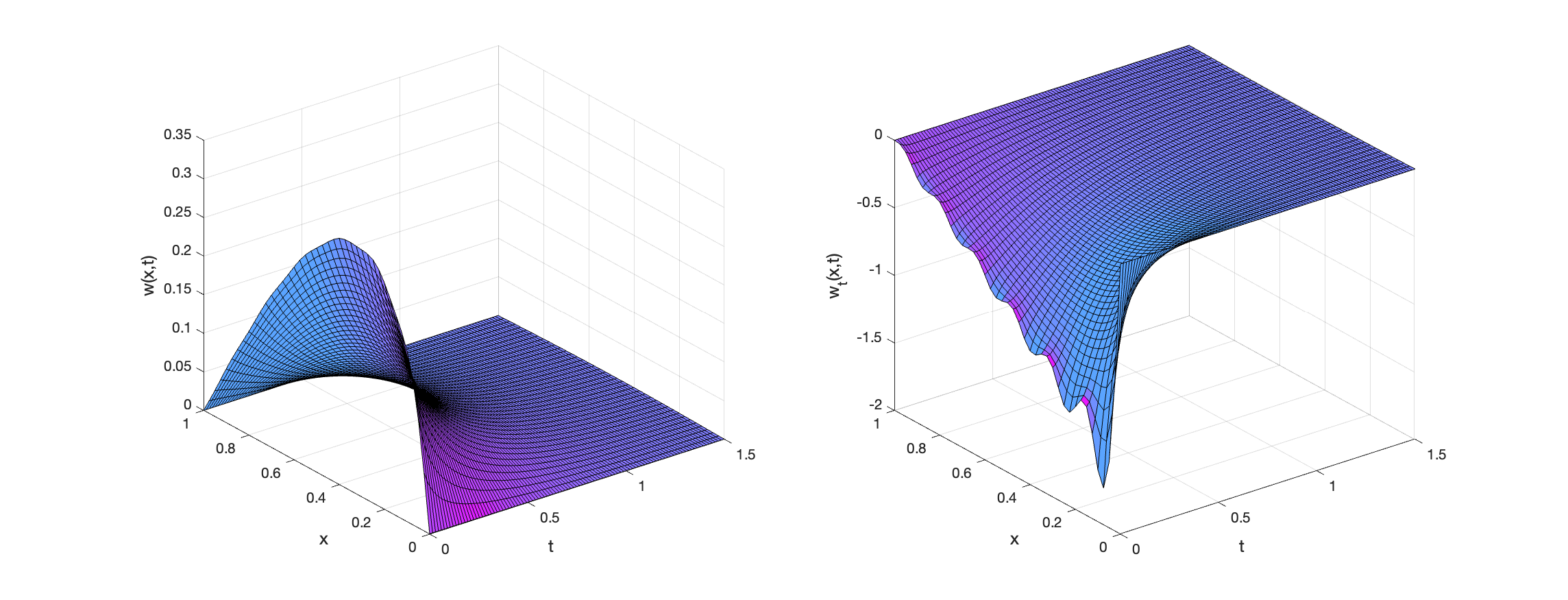}
\end{flushleft}
\caption{Functions $w$ and $w_t$ for $L=12$}
\label{figW12}
\end{figure}

\begin{figure}
\begin{flushleft}
\includegraphics[height=8.9cm,width=16cm]{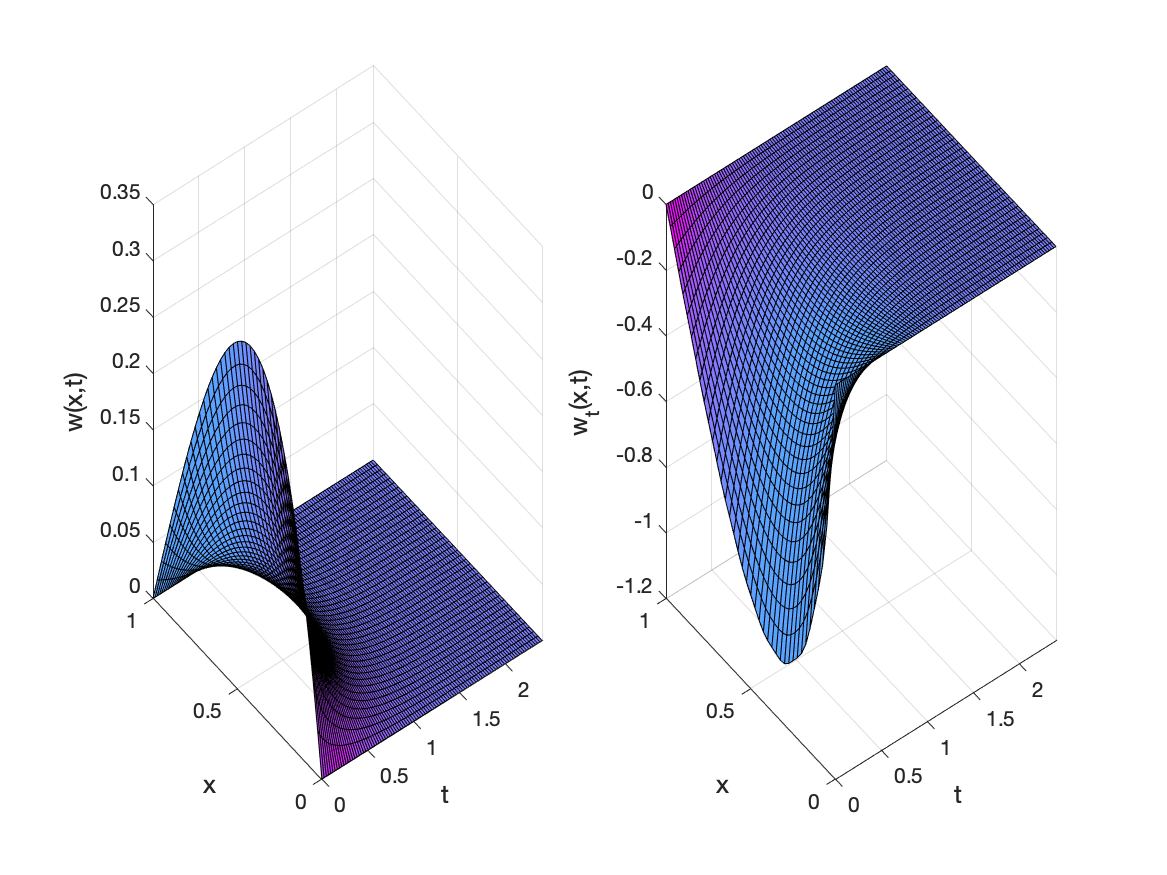}
\end{flushleft}
\caption{Functions $w$ and $w_t$ for $L=16$}
\label{figW16}
\end{figure}

Numerical verification of (\ref{a7}) and (\ref{a8}) was performed
by means of  test problems, where $\Omega=(0,1)$,
$$
\phi_i(x)=\sqrt{2}\sin (i\pi x), \quad \lambda_i=i^2\pi^2, \qquad 
f(x)=\sum\limits^L_{i=1}\frac{1}{i^m}\sin(i\pi x).
$$
The corresponding exact solution $w(x,t)$ and its derivative $w_t(x,t)$  are 
depicted in Fig.~\ref{figW12}
(for $m=1$ and $L=12$) and  in Fig.~\ref{figW16}
(for $m=2$ and $L=16$).

\begin{figure}
\centerline{\includegraphics[width=16cm,height=8cm]{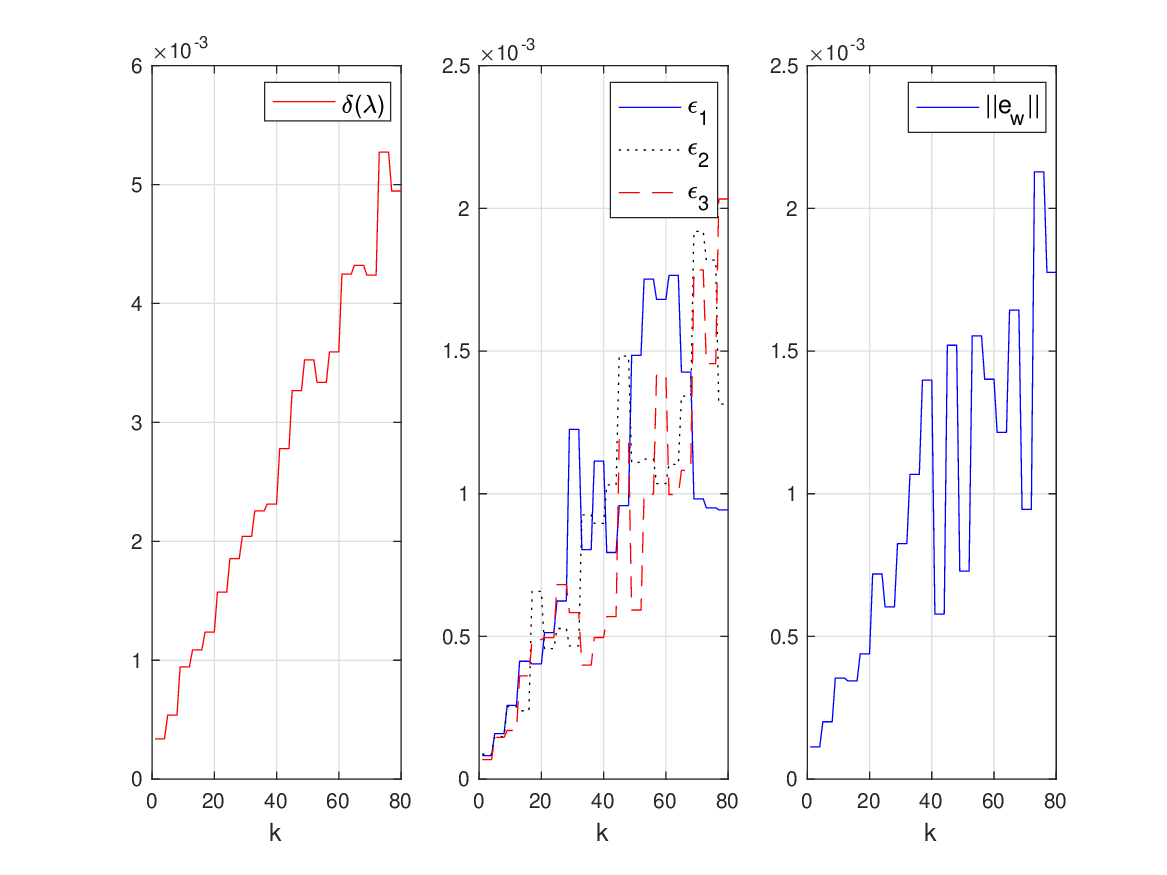}}
\caption{Disturbances of exact eigenvalues and eigenfunctions and the overall 
error for $L=N=12$.}
\label{figDEL1212}
\end{figure}

\begin{figure}
\centerline{\includegraphics[height=8cm,width=6.9cm]{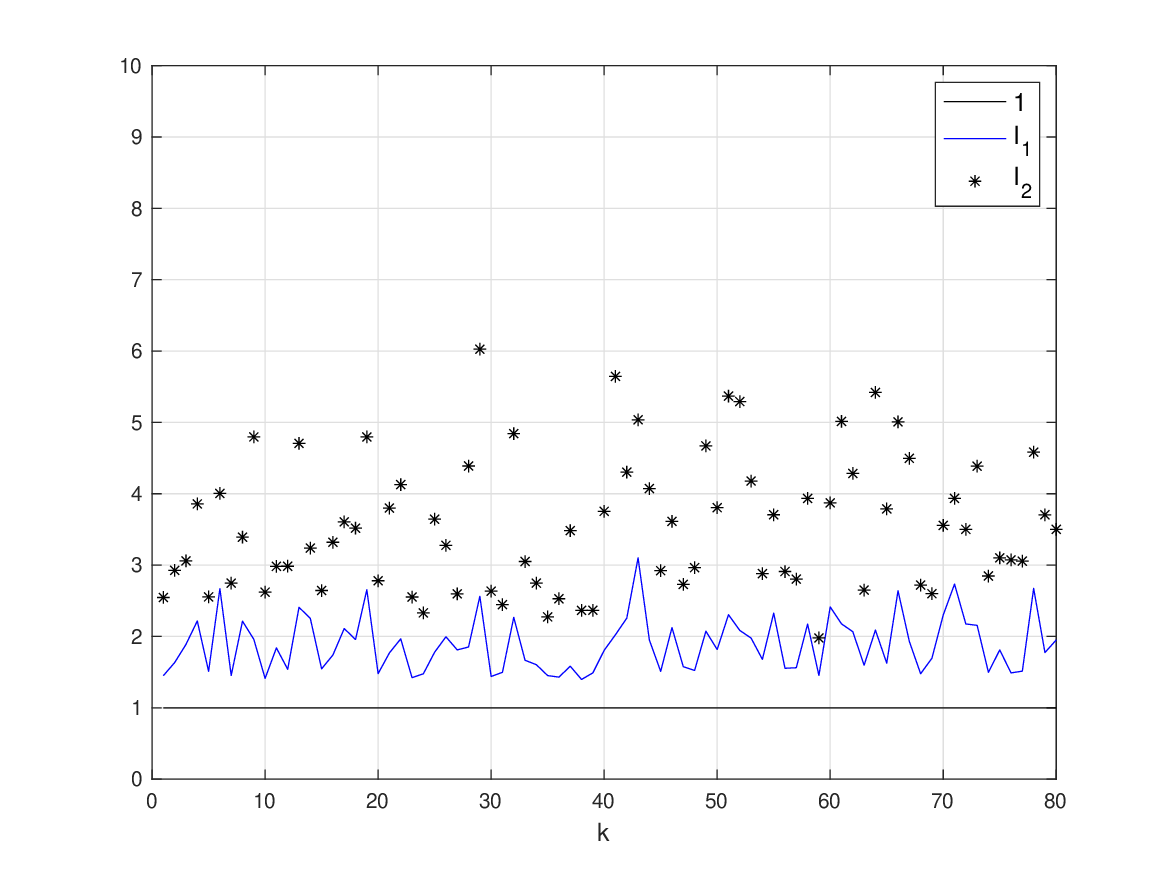}\!\!\!
\includegraphics[height=8cm,width=6.9cm]{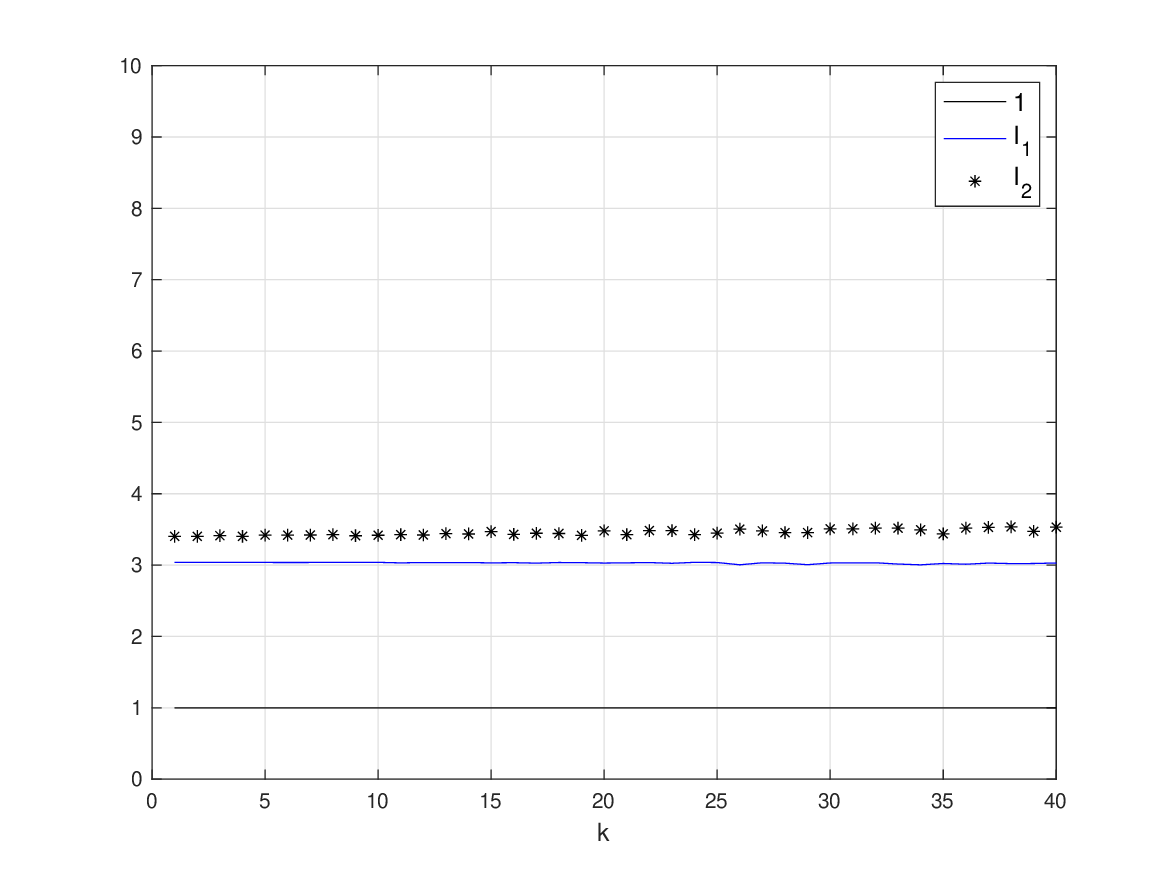}}
\caption{Indexes of the Majorant }
\label{figIND1}
\end{figure}

Approximation $\wt w$
is defined in accordance with (\ref{a4}) and $S_N$ is computed by the relation 
(\ref{SN}). The numbers  $\theta_i$ are obtained by random disturbances 
of $\lambda_i$ and the amplitude of these disturbances is increasing with the 
number $k$ of a test ($k=1,2,...,n$). Analogously, approximate eigenfunctions 
$\psi_i$ are obtained
by certain disturbances of $\phi_i$. 
Disturbances are characterised by the 
quantities
\ben
\label{del}
\delta=\frac{1}{L}\sum\limits^L_{i=1}\frac{|\lambda_i-\theta_i|}{\lambda_i}
\quad\text{ and}\quad \epsilon_i=\|\phi_i-\psi_i\|_{L_2},\quad i=1,2,...,L.
\een
The left picture in Fig. \ref{figDEL1212} shows errors encompassed
in $\theta_i$ (for $L=12$)  and the central picture depicts errors associated
with the eigenfunctions $\psi_1$,  $\psi_2$, and $\psi_3$.
In this example, $\wt w$ is defined by (\ref{a4}), where $N=12$.
Therefore, the overall error (right picture in Fig.  \ref{figDEL1212}) is
generated exclusively by inaccuracies in eigenvalues and eigenfunctions.
In Fig. \ref{figIND1} we depict quantities that characterise efficiency
of error estimates (\ref{a7}) and (\ref{a8}).
Here $I_1$ is the ratio between $M_\oplus(\wt w;\y)$ and
$\Normt{\nabla_{xt}\e_w}$, while 
$I_2$ is the square root of the analogous ratio between the right
and left parts of (\ref{a8}) (for $\alpha=0.5$). Mean values of $I_1$ and $I_2$
in this test are presented in the first line of Tab. \ref{tabA}.

  There is another factor, which may affect the accuracy of $\wt w$.
  It is the value of $N$. 
  If the approximation (\ref{a4})  contains an insufficient number of terms, 
then $\wt w$ may
  contain a significant error even if the eigenvalues and eigenfunctions are 
defined exactly.
  This situation is seen in the test series shown in the fourth line of Tab. 
\ref{tabA}
  and the
  corresponding indexes are depicted in Fig. \ref{figIND1} (right). 
   Here,  the main part of the error arises because $N<L$ and, therefore, the 
error
   is insensitive to changes in determining the eigenvalues.

\begin{table}[h!]
\begin{center}
\begin{tabular}{||c|c|c|c|c|c|c|c|c||}
\hline
$n$\;&\;$m$\;&\;$L$\;\;&\;\;$N$\; \;&\;\;$\alpha$\; \;&\;\;$I_1$\;\; 
&\;\;$I_2$\;\;
&\;\;$\delta_{max}\;\;$&\;\;$\epsilon_{max}$\;\;\\[2pt]
\hline
80& 1 &  12 &  12   & 0.5 & 1.888 & 3.501  &0.003&0.015\\
80& 2  &   16&  12   & 0.5 & 2.215  & 3.129   & 0.002  & 0.015\\
20 & 1.5 &   8 &  4   & 0.5 &3.014  & 3.476 & 0.003&0.010\\
40 &1.5& 8 &  4   & 0.5  &3.030  & 3.457 &0.004&0.015 \\
20 &2& 8 &  8  & 0.5     &1.827  &  3.355  &0.007&0.048\\
20 &1& 10 &  10   & 0.3 &2.029  &  3.351 &0.012&0.048 \\
20 &1& 10 &  8   & 0.3.   &4.782  &  4.411 &0.004&0.024 \\
20 &1.5& 12 &  12   & 0.3 &1.828  &  2.566 &0.004&0.024 \\
   \hline
\end{tabular}
\end{center}
\vskip3pt
\caption{Averaged indexes  and other data for different series of tests.}
\label{tabA}
\end{table}
Also, Tab. \ref{tabA} exposes results related to  other experiments.
The first column shows the amount of tests in each series, parameter $m$
is used in the definition of $f$, and parameter $\alpha$ is  used in the 
estimate (\ref{a8}).
The last two columns show maximal disturbances of  eigenvalues and 
eigenfunctions
(see (\ref{del})). Different efficiency indexes
$I_1$ and $I_2$ show the quality of estimates. If an estimate is sharp, then 
such
 type index is close to one. 
Tab.  \ref{tabA} shows that estimates admit certain overestimation (typically 
from 2 to 5 times).
As noted earlier, the estimates were obtained using the simplest method.
However, they give a correct presentation on  the actual values of errors.
Therefore, we may await that using more sophisticated techniques will give
in much sharper results.


\end{document}